\documentclass[12pt]{article}
\usepackage{amsmath,amsthm,latexsym,amssymb,amscd,amsfonts,dsfont,txfonts,times}
\usepackage{cases,color,graphicx,epsfig,mathrsfs}
\usepackage[numbers,sort&compress]{natbib}
\allowdisplaybreaks
\usepackage{mathtools}
\usepackage{amsmath}
\usepackage{amsthm}
\usepackage{amsfonts}
\usepackage{bm}
\usepackage{mathrsfs}
\usepackage[all]{xy}
\usepackage{color}
\usepackage[rgb]{xcolor}
\usepackage{float}
\usepackage{subfig}
\usepackage{bookmark}
\usepackage{multirow}
\usepackage{diagbox}
\usepackage{booktabs}
\usepackage[font=footnotesize]{caption}
\numberwithin{equation}{section}
\newtheorem{theorem}{Theorem}[section]
\newtheorem{lemma}[theorem]{Lemma}
\newtheorem{remark}{Remark}[section]

\newtheorem{example}{Example}[section]
\newcommand{\commentOne}[1]{#1}
\newcommand{\commentTwo}[1]{#1}
\newcommand{\ctone}[1]{#1}
\newcommand{\commentone}[1]{#1}
\newcommand{\commenttwo}[1]{#1}
\newcommand{\commentthree}[1]{#1}
\title{\Large{\bf Effective numerical treatment of sub-diffusion equation with non-smooth solution}}
\author{
   \small Zongze Yang$^a$,\quad
   \small  Jungang Wang$^a$\,\footnote{Corresponding author. E-mail address: modiker@163.com}\ ,\quad
   \small Yan Li$^b$,\quad
    and\  \small  Yufeng Nie$^a$\\
   {\small\it $^a$Department of Applied Mathematics, Northwestern Polytechnical University, P. R. China}  \\
   {\small\it $^b$Department of Mathematics, City University of Hong Kong, Hong Kong, P. R. China}\\
}
\date{ }
\begin{document}
\maketitle
\begin{abstract}
In this paper we investigate a sub-diffusion equation for simulating the
anomalous diffusion phenomenon in real physical environment. Based on an
equivalent transformation of the original sub-diffusion equation followed by
the use of a smooth operator, \commentTwo{we devise a high-order numerical scheme by
combining the Nystr\"om method in temporal direction with the compact finite
difference method and the spectral method in spatial direction.} The distinct
advantage of this approach in comparison with most current methods is its high
convergence rate even though the solution of the anomalous sub-diffusion
equation usually has lower regularity on the starting point. The effectiveness
and efficiency of our proposed method are verified by several numerical
experiments.
\end{abstract}
\section{Introduction}
Fractional calculus is an area having a long history, which is believed to have
stemmed from a question about the meaning of notation
$d^\frac{1}{2}y/dx^\frac{1}{2}$ raised in the year 1695 by Marquis de
L'H\^opital to Gottfried Wihelm Leibniz. During the past three decades, this
subject has gained considerable popularity due mainly to its powerful
applications in numerous seemingly diverse and widespread fields of science and
engineering, such as materials and mechanics, signal processing, anomalous
diffusion, biological systems, finance, 
and etc.(see \cite{Zaslavsky2002, Yuste2004, Metzler2000, Cartea2007, Magin2010, Gorenflo2001}).
At present there have been many papers presenting fractional calculus models
for kinetics of natural anomalous processes in complex systems.
\commentTwo{These models always maintain the long-memory and non-local
properties of the corresponding dynamics. Because of these properties, it is
still not easy to find the exact or numerical solutions of these equations,
though researchers have developed many methods to approach this goal.}
Special interest has been paid to the anomalous diffusion processes, which include
super-slow diffusion (or sub-diffusion) and super-fast diffusion (or super-diffusion).
Among those models, anomalous sub-diffusion equations are important due to its
application in simulating real physical sub-diffusion phenomena. 
\commentthree{%
The model is always written as
  \begin{equation*}
    \left\{
    \begin{aligned}
      & u_t(x,t) = \prescript{RL}{0}{\mathcal D}^{1-\gamma}_t Lu(x,t) + \bar f(x, t), &0 < x \le X, \, 0 < t \le T, \\
      &u(x,0) = \varphi(x), &0 \le x \le X, \\
      &u(0,t) = \psi_1(t), \quad u(X,t) = \psi_2(t),&0 \le t \le T, \\
    \end{aligned}
    \right.
  \end{equation*}
  where $ 0 < \gamma< 1$, $Lu(x,t) = K_{\gamma} u_{xx}(x,t)$ and
$\prescript{RL}{0}{\mathcal D}^{1-\gamma}_t $
denotes the Riemann-Liouville fractional derivative of order $1-\gamma$,
\begin{equation}
  \prescript{RL}{0}{\mathcal D}^{1-\gamma}_t u(x, t) =
  \frac{1}{\Gamma(\gamma)} \frac{\partial}{\partial t}\int_0^t (t-\eta)^{\gamma-1}u(x,\eta)\,d\eta.
\end{equation}
\commentTwo{Assume $\bar f= \prescript{RL}{0}{\mathcal D}^{1-\gamma}_t f$,
then we can rewrite the original equation as below\cite{Metzler2000}}:
  \begin{equation}\label{eq:orig}
    \left\{
    \begin{aligned}
      & u_t(x,t) = \prescript{RL}{0}{\mathcal D}^{1-\gamma}_t \big(Lu(x,t) + f(x, t)\big), &0 < x \le X, \, 0 < t \le T, \\
      &u(x,0) = \varphi(x), &0 \le x \le X, \\
      &u(0,t) = \psi_1(t), \quad u(X,t) = \psi_2(t),&0 \le t \le T.
    \end{aligned}
    \right.
  \end{equation}
In this paper we only consider anomalous sub-diffusion problem in form~\eqref{eq:orig}.}
In some references, equation~\eqref{eq:orig} is called the time Riemann-Liouville type
sub-diffusion equation. \commentthree{Some researchers use the following models instead of
equation~\eqref{eq:orig} :
\begin{equation}
 \prescript{RL}{0}{\mathcal D}^{\gamma}_t (u(x,t)-u(x,0)) =  Lu(x,t) + f(x, t),
\end{equation}
or 
\begin{equation}
 \prescript{C}{0}{\mathcal D}^{\gamma}_t u(x,t)=  Lu(x,t) + f(x, t),
\end{equation}
where $\prescript{C}{0}{\mathcal D}^{\gamma}_t$ denotes the Caputo fractional derivative of order $\gamma$,
\begin{equation}
  \prescript{C}{0}{\mathcal D}^{\gamma}_t u(x, t) =
  \frac{1}{\Gamma(1-\gamma)} \int_0^t (t-\eta)^{-\gamma}\frac{\partial}{\partial \eta}u(x,\eta)\,d\eta.
\end{equation}}
In fact, these models are equivalent.
Many numerical methods have been developed to solve anomalous sub-diffusion
equations.  In 2005, Yuste and Scedo\cite{Yuste2005} proposed an explicit FTCS
scheme, which combined the forward time centered space (FTCS) method with the
Gr\"unwald-Letnikov discretization of the Riemann-Liouville derivative. And a
new von Neumann-type method is applied to analysis the stability in the paper.
Zhuang et al.\cite{Zhuang2008} presented an implicit numerical method as well as 
two techniques which are used to improve the order of convergence. 
The stability and convergence
analysis for the implicit numerical method are given by using an energy method. 
\commentTwo{Combining the L1 discretization for time-fractional part and
fourth-order accuracy compact approximation for space derivative, 
a compact finite difference scheme is established by 
Gao and Sun\cite{Gao2011}.} 
Furthermore, Gao et al.\cite{Gao2015} offered a scheme with global second-order
numerical accuracy in time independent of the fractional derivative exponent.
\commentTwo{Apart from finite difference methods, part of researchers have
  investigated Galerkin methods, including finite element methods and spectral
methods.}
Zeng et al.\cite{Zeng2013} adopted linear
multistep method and finite element method to approach the time-fractional
sub-diffusion equation and got two unconditionally stable schemes.
\commentthree{Mustapha developed a discontinuous Petrov-Galerkin 
  method\cite{Mustapha2014} and a time-stepping $hp$-versions discontinuous Galerkin
  method\cite{Mustapha2015} for
time fractional partial differential equations.}
\commenttwo{%
  In \cite{Dehghan2017, Dehghan2016}, Dehghan et al.\ developed spectral element
  method in spatial and finite difference method in temporal for nonlinear fractional partial 
  differential equations 
  and sub-diffusion equations, and gave the corresponding theoretical analysis. 
  They also studied the homotopy analysis method and the dual reciprocity boundary integral method for
  fractional partial differential equations\cite{Dehghan2015a, Dehghan2010}.
  In addition, the authors presented two high-order methods for multi-term
  time-fractional diffusion equation\cite{Dehghan2015b}.}
\commentTwo{For multi-term time-fractional diffusion equation, with the benefit
  of spectral method, Zheng et al.\cite{Zheng2015} gained a valuable high-order
  scheme, which possessing high efficiency and exponential decay in both time
and space directions.}
There have been a great deal of researches on anomalous sub-diffusion
equations, however, as noticed in~\cite{Gao2015},
when the solution is not smooth enough at $t=0$, \commentthree{the convergence rate will
be lower than expectation}.
To overcome this shortcoming, we adopt some techniques similar to those  used
to deal with weakly singularity Volterra integral equation\cite{Baratella2004, Monegato1997}.
\commentOne{By these techniques, we can obtain better numerical results even though the
  solution has weak regularity.  The effectiveness of our algorithm
can be seen in the numerical examples.}
The outline of this paper is arranged as follows. In Section~\ref{sec:ts}, we
give an equivalent form of equation~\eqref{eq:orig} by equivalent
transformation and smoothing operator, which can improve the regularity of the solution.
\commentTwo{Section~\ref{sec:scheme} contains a semi-discrete scheme given by
discretizing the equivalent equation with Nystr\"om method. 
With different method discretizing spatial variables, two fully-discrete
schemes are presented in \commentone{Section}~\ref{sec:full}.
To demonstrate the efficiency and effectiveness of the proposed scheme,
we perform some numerical examples in Section~\ref{sec:example}. 
And in the last section conclusions as well as some remarks are given.}
\section{Equivalent transformation and smoothing method}\label{sec:ts}
In paper \cite{Baratella2004}, the authors proposed a simple smooth
transformation  for Volterra integral equations, which can improve the
regularity of the solution and can be used to construct high-order convergence 
methods.  \commentthree{In this section, we further introduce a smoothing method
to transform the original fractional differential
equation~\eqref{eq:orig} into an equivalent form}.
In order to use the smoothing method, the equation \eqref{eq:orig} is transformed into
an integral equation by integrating both sides:
\commenttwo{
\begin{equation*}
  \int_0^s u_t(x, t)\,dt = \int_0^s \frac{1}{\Gamma(\gamma)}\frac{\partial}{\partial t}\int_0^{t} (t-\eta)^{\gamma-1}
    \big(Lu(x, \eta) + f(x, \eta)\big)\,d\eta \,dt.
\end{equation*}
Then we have
\begin{equation*}
  u(x, s) = u(x, 0) + \frac{1}{\Gamma(\gamma)}\int_0^{s} (s-\eta)^{\gamma-1}
  \big(Lu(x, \eta) + f(x, \eta)\big)\,d\eta,
\end{equation*}
i.e.}
\begin{equation}\label{eq:int}
  u(x, t) = u(x, 0) + \frac{1}{\Gamma(\gamma)}\int_0^{t} (t-\eta)^{\gamma-1}
    \big(Lu(x, \eta) + f(x, \eta)\big)\,d\eta.
\end{equation}
The last term in the right hand of \eqref{eq:int} has a similar form with an integral term in the Volterra integral equation.
Following \cite{Baratella2004}, we introduce the smooth operator
\begin{equation}
  \lambda(t) = (b - a)^{1-q}(t-a)^q+a, \quad q \in \{1, 2, \dots,n,\dots\},
\end{equation}
which maps $[a, b]$ into $[a, b]$ \commentTwo{where $a$ is $0$ and $b$ is $T$}.
\commentTwo{Here, we use $a$, $b$ as end points to state the generality of the transformation.}
Let  $\alpha = 1-\gamma$ and change the variables in \eqref{eq:int} 
by setting $\eta = \lambda(s)$,
$t = \lambda(t')$. Replacing  $t'$ by $t$, we then get
\begin{equation}\label{eq:intab}
  \begin{aligned}
    u(x, \lambda(t)) =  u(x, \lambda(a)) +
    \frac{1}{\Gamma(\gamma)}\int_a^{t}\big(\lambda(t)-\lambda(s)\big)^{-\alpha}
      G(x, s)\lambda'(s)\, ds, \\
  \end{aligned}
\end{equation}
where $G(x, s) =Lu\big(x,\lambda(s)\big)+f\big(x,\lambda(s)\big)$.
To transform the kernel of \eqref{eq:intab} with form $(t-s)^{-\alpha}$, 
we denote\cite{Baratella2004}
\begin{equation}\label{eq:da}
  \delta_{\alpha}(t,s) =  \left\{
    \begin{aligned}
      &\left( \frac{(t-a)^q-(s-a)^q}{t-s}\right)^{-\alpha}, &t &\ne s,\\
      &\big(q(s-a)^{q-1}\big)^{-\alpha}, &t &= s.
    \end{aligned}
    \right.
\end{equation}
The equation \eqref{eq:da} implies
\begin{equation}
  \begin{aligned}
      \big(\lambda(t)-\lambda(s)\big)^{-\alpha}
      &= \big((b-a)^{1-q}\big)^{-\alpha}\delta_{\alpha}(t, s)(t-s)^{-\alpha}.
  \end{aligned}
\end{equation}
Multiplying both sides of \eqref{eq:intab} by $\lambda'(t)$, we can now rewrite \eqref{eq:intab} as
\begin{equation}\label{eq:intab2}
  \begin{aligned}
  \lambda'(t)u(x, \lambda(t)) = \lambda'(t)u(x, \lambda(a))
  + \frac{1}{\Gamma(\gamma)}\int_a^{t}(t-s)^{-\alpha}K_\alpha(t,s)
  G(x,s)\lambda'(s)\,ds,
  \end{aligned}
\end{equation}
where
\begin{equation}
  K_{\alpha}(t, s) = \big((b-a)^{1-q}\big)^{-\alpha}\lambda'(t)\delta_{\alpha}(t, s).
\end{equation}
In order to use the Nystr\"om method in spatial direction, we introduce another
transformation $\mu(t) = \frac{b-a}{2}t + \frac{b+a}{2}$
and denote
\begin{equation}
  \left\{
  \begin{aligned}
    v(x, t) &= \lambda'\big(\mu(t)\big)u\big(x, \lambda\big(\mu(t)\big)\big), \\
    g(x, t) &= \lambda'\big(\mu(t)\big)f\big(x, \lambda\big(\mu(t)\big)\big), \\
    h(x, t) &= \lambda'\big(\mu(t)\big)u\big(x, \lambda\big(a\big)\big).
  \end{aligned}
  \right.
\end{equation}
By setting $t=\mu(t')$, $s=\mu(s')$ in \eqref{eq:intab2} and replacing
$t'$ by $t$, $s'$ by $s$, we get
\begin{equation}\label{eq:smooth}
  v(x, t) = h(x, t) + \frac{1}{\Gamma(\gamma)}\int_{-1}^{t}
      (t-s)^{-\alpha}H(t, s)\big(Lv(x,s)+g(x,s)\big)\,ds,
\end{equation}
where
\begin{equation}
  H(t, s) = \left(\frac{b-a}{2}\right)^{1-\alpha} K_\alpha\big(\mu(t), \mu(s)\big).
\end{equation}
 After two times of transformations, a new equation~\eqref{eq:smooth} for the problem~\eqref{eq:orig}
 is derived. 
\commenttwo{%
For equation \eqref{eq:smooth}, the kernel is $\hat H(t, s) = (t-s)^{-\alpha}H(t,s)$.
According to the expression, we know $H(t,s)$ is continuous and 
$\hat H(t, s) \in L^1([-1,1]\times[-1,1])$ is a weakly singular kernel. %
}
Now we consider the smoothness of the solution of the new equation \eqref{eq:smooth}.
\commentOne{%
Suppose that $u(x,t) \to O(t^\beta)$ ($\beta>0$) as $t \to 0$. 
Then $\lambda'(s) u(x, \lambda(s)) \to O(s^{q\beta+q-1})$ as $s \to 0$.
As we set $v(x,t) = \lambda'\big(\mu(t)\big)u\big(x,\lambda(\mu(t))\big)$ 
and $\mu(t)$ is the first order polynomial which do not change the smoothness.
Then $v(x, t) \to O(t^{q\beta+q-1})$ as $t \to 0$. Obviously, $q\beta +q -1 \ge \beta$, i.e.\ 
$v(x, t)$ is smoother than $u(x, t)$ when $q > 1$. Because  
the new solution has higher regularity, traditional methods can be applied to
solve this equation without loss much accuracy.}
As a conclusion, we illustrate the smoothing process as below
\begin{displaymath}
  \xymatrixcolsep{8pc}
  \xymatrix{%
    u(x,t) \ar[r]^-{\lambda(s)=(b-a)^{1-q}(s-a)^q+a} &
    u(x,\lambda(s)) \ar[d]^{\times \lambda'(s)} \\
    v(x,r)&\lambda'(s)u(x,\lambda(s))
    \ar[l]_-{s=\mu(r)}^-{\mu(r)=\frac{b-a}{2}r+\frac{b+a}{2}} }
\end{displaymath}
\commentTwo{In the next section, we will describe more details of how to
discretize the above equation.}
\section{The semi-discrete approximation}\label{sec:scheme}
In the previous section we have obtained \commentone{an equivalent equation~\eqref{eq:smooth}}.
Here, we present a semi-discrete scheme by using the Nystr\"om method.
Choose $N+1$ distinct points $\tau_n, n = 0, \cdots, N$ in the interval
$[-1, 1]$ corresponding to $t_n, n = 0, \cdots, N$ in $[a, b]$ with
$t_n = \lambda\big(\mu(\tau_n)\big)$. And collocate the equation at the nodes $\{\tau_n\}^{N}_{n=0}$
\begin{equation}
  v(x, \tau_n) = h(x, \tau_n) + \frac{1}{\Gamma(\gamma)}\int_{-1}^{\tau_n}
      (\tau_n-s)^{-\alpha}H(\tau_n, s)\big(Lv(x,s)+g(x,s)\big)ds.
\end{equation}
Then replace $H(\tau_n, s)\big(Lv(x,s)+g(x,s)\big)$ by the corresponding
Lagrange interpolation polynomials associated with $\{\tau_n\}^{N}_{n=0}$
\begin{equation}
  \begin{aligned}
  v(x, \tau_n) = h(x,\tau_n) + \frac{1}{\Gamma(\gamma)}
    \sum_{j=0}^{N}w_{n, j}H(\tau_n, \tau_j)(Lv(x, \tau_j) + g(x, \tau_j))
  \end{aligned}
\end{equation}
where
$$w_{n, j} = \int_{-1}^{\tau_n} (\tau_n-s)^{-\alpha}I_{N,j}(s)\,ds,$$
and
$I_{N,j}(s)$ are the Lagrange interpolation polynomials on points $\{\tau_n\}^{N}_{n=0}$.
For the computation of the coefficients of $w_{n, j}$, we use the Jacobi-Gauss
quadrature as described in \cite{Shen2011}.
\begin{remark}
  \commenttwo{Because of the singularity of the solution, the point $\tau = -1$ should not be chosen as an element of $\{\tau_n\}$.}
\end{remark}
\begin{remark}
  By choosing different $\{\tau_n\}_{n=0}^{n=N}$, we get different approximation polynomials with different accuracies.
  This will have effect on the accuracy of the solution.
\end{remark}
Now, we obtain a semi-discrete scheme for equation \eqref{eq:orig}
\begin{equation}\label{eq:semischeme}
  \begin{aligned}
  v(x, \tau_n) = h(x,\tau_n) +
  \sum_{j=0}^{N}r_{nj}(Lv(x, \tau_j) + g(x, \tau_j)), \quad n=0,1,\cdots, N,
  \end{aligned}
\end{equation}
where $r_{nj} = w_{n, j}H(\tau_n, \tau_j)/\Gamma(\gamma) $. Once $v(x, \tau)$ is obtained, the solution
of the original equation is given by
\begin{equation}\label{eq:trans}
  u(x, t) = \frac{v(x, \tau)}{\lambda'\big(\mu(\tau)\big)}, \quad
  t = \lambda\big(\mu(\tau)\big).
\end{equation}
\section{The fully-discrete approximation}\label{sec:full}
Because of the high convergence rate in temporal direction of the scheme,
we need a high-order method in spatial direction for competence.
In this section, we perform two different discrete methods on the spatial variable.
\subsection{Discrete spatial variable by compact difference operator} \label{sec:improve}
In this subsection, we use spatial compact approximation in spatial direction.
Let $x_k = k\Delta x,(k = 0, 1, \cdots, M)$ with step $\Delta x = X/M$.
Like~\cite{Gao2015}, denote a average operator
\begin{equation}
  \mathcal{A}u_i = \left\{
    \begin{aligned}
      &(I+\frac{\Delta x^2}{12}\delta^2_x)u_i, & 1 \le i \le M-1, \\
      &u_i,& i =0 \  \text{or} \ i = M.
    \end{aligned}
    \right.
\end{equation}
where $\delta_x^2$ is the centered difference operator.
Perform $\mathcal{A}$ on both sides of \eqref{eq:semischeme} at $\{x_k\}_{k=1}^{M-1}$
\commentone{%
\begin{equation}
  \begin{aligned}
    \mathcal{A}v(x_k, \tau_n) = \mathcal{A}h(x_k,\tau_n) +
    \sum^{N}_{j=0}r_{nj} \big(K_{\gamma}\mathcal{A}v_{xx}(x_k,\tau_j)
    +\mathcal{A}g(x_k,\tau_j)\big). \\
  \end{aligned}
\end{equation}%
}%
The next thing is to approximation $\mathcal{A}v_{xx}(x_k,s)$.
In order to obtain the spatial compact scheme, we need the following lemma, which suggests that
$\delta_x^2v(x_k, s)$ is a good approximation to $\mathcal{A}v_{xx}(x_k,s)$.
\begin{lemma}[\cite{Gao2015}]\label{lemma:c}
  Let function $g(x) \in C^6[0, X]$ and
  $\xi(\eta)=(1-\eta)^3\big(5 - 3(1-\eta)^2\big)$. Then
  \begin{equation}
    \begin{aligned}
    \mathcal{A}g''(x_i) = \delta^2_xg(x_i) + \frac{\Delta x^4}{360}\int_0^1\big(
      g^{(6)}(x_i-\eta \Delta x) + g^{(6)}(x_i + \eta\Delta x)\big)\xi(\eta)\,d\eta.
    \end{aligned}
  \end{equation}
\end{lemma}
By Lemma~\ref{lemma:c}, we have
\begin{equation}
  \begin{aligned}
    \mathcal{A}v(x_k, \tau_n) = \mathcal{A}h(x_k, \tau_n) +
    r_{nj}(K_{\gamma}\delta^2_xv(x_k, \tau_j) + \mathcal{A}g(x_k, \tau_j)) + O(\Delta x^4).
  \end{aligned}
\end{equation}
Drop down the high-order term, hence
\begin{equation}\label{eq:compact}
  \begin{aligned}
    \mathcal{A}v(x_k, \tau_n) = \mathcal{A}h(x_k, \tau_n) +
    r_{nj}(K_{\gamma}\delta^2_xv(x_k, \tau_j) + \mathcal{A}g(x_k, \tau_j)).
  \end{aligned}
\end{equation}
Now, we have established the fully-discrete scheme by discrete spatial variables with compact difference operator.
To introduce the matrix form of the last scheme, we denote $v_k^n = v(x_k, \tau_n)$, $h_k^n = h(x_k, \tau_n)$, and $g_k^j = g(x_k, \tau_j)$.
\commentTwo{%
In addition, let $\bm{V}$, $\bm{V}_0$, $\bm{W}$ and $\bm{G}$ be matrix with  
    $ (\bm{V})_{k, n+1} = v_k^n$, $(\bm{V}_0)_{k, n+1} = h_k^n$, $(\bm{W})_{n+1,j+1} =r_{nj}$ and
    $(\bm{G})_{k, n+1}= g_k^n$\,
    ($k = 1,\cdots, M-1$, \,$j = 0, \cdots, N$, \,$n = 0, \cdots, N$).
    Then the matrix form of scheme \eqref{eq:compact} can be written
    as below:
\begin{equation}\label{eq:matrixc}
  \begin{aligned}
    \bm{V}+\frac{1}{12}(\bm{DV}+\bm{B}_v) =\,&\bm{V}_0 + \frac{1}{12}(\bm{D}\bm{V}_0+\bm{B}_h) 
  +\frac{K_{\gamma}}{\Delta x^2}(\bm{D}\bm{V}+\bm{B}_v)\bm{W}^T \\ &+
  \big(\bm{G}+\frac{1}{12}(\bm{D}\bm{G}+\bm{B}_g)\big) \bm{W}^T,
\end{aligned}
\end{equation}%
where the matrix $\bm{D}$ equals to $\text{tridiag}(1, -2, 1)$,
\commentthree{ and $\bm{B}_h$, $\bm{B}_g$ are similar to $\bm{B}_v$} which 
is defined as follows:
\begin{equation*}
  \begin{aligned}
    &\bm{B}_{v,(M-1)\times(N+1)} =  \left(
    \begin{array}{ccccc}
      v_0^0  & v_0^1 & \cdots & v_0^{N-1} & v_0^N \\
            0        &      0        & \cdots &       0         &      0        \\
        \vdots       &    \vdots     & \ddots &   \vdots        &    \vdots     \\
            0        &      0        & \cdots &       0         &      0        \\
      v_M^0  & v_M^1 & \cdots & v_M^{N-1} & v_M^N
    \end{array}
    \right).  \\
  \end{aligned}
\end{equation*}
.}
For simplification, we recombine the terms in equation~\eqref{eq:matrixc} as
\begin{equation}\label{eq:matrixcc}
  \bm{T}\bm{V} - \bm{A}\bm{V}\bm{W}^T = \bm{S},
\end{equation}
where $ \bm{T} = \bm{I} + \frac{1}{12}\bm{D}$,
$ \bm{A} = \frac{K_\gamma}{\Delta x^2}\bm{D}$ and
\commentthree{%
\begin{equation}
\bm{S} = -\frac{1}{12}\bm{B}_v + \bm{V}_0 + \frac{1}{12}(\bm{D}\bm{V}_0+\bm{B}_h)
+ \big(\bm{G}+\frac{1}{12}(\bm{D}\bm{G}+\bm{B}_g)\big) \bm{W}^T.
\end{equation}%
}%
Combining \eqref{eq:trans} and \eqref{eq:matrixcc}, we can obtain the
solution of the original equation.
\subsection{Spatial discretization with spectral method}\label{sec:spectral}
In order to compare the result with some existing algorithms, we use the space basis
function, presented by Zheng et al.\ in \cite{Zheng2015}, to solve equation~\eqref{eq:semischeme}.
Here, we only consider zero boundary condition, i.e., $u(0, t) = 0$, $u(X,t) = 0.$
Firstly, let us state the variational problem of equations~\eqref{eq:semischeme}.
We recall the semi-discrete scheme~\eqref{eq:semischeme}
\begin{equation}
  \begin{aligned}
    v(x, \tau_n) = h(x,\tau_n) + \sum_{j=0}^{N}r_{nj}(Lv(x,\tau_j)  + g(x, \tau_j))\quad
    n=0,1,\cdots, N.
  \end{aligned}
\end{equation}
Set $v^n(x) = v(x,\tau_n)$, $h^n(x) = h(x,\tau_n)$, $g^n(x) = g(x,\tau_n)$,
$\Lambda =[0,X]$. Then the variational problem
of \eqref{eq:semischeme} is :
Find $v^n(x) \in H_0^1(\Lambda)$, $(n = 0, 1, \cdots, N)$ for
$\forall \phi(x) \in H_0^1(\Lambda)$, such that
\begin{equation}\label{eq:var}
  \begin{aligned}
    \big(v^n(x),\phi(x)\big) + &\ K_\gamma\sum_{j=0}^{N}r_{nj}\big(v^j_x(x),\phi_x(x)\big)  \\
    &=\big(h^n(x),\phi(x)\big) +  \sum_{j=0}^{N}r_{nj}\big(g^j(x),\phi(x)\big).
  \end{aligned}
\end{equation}
Next, we construct a spectral scheme for the above variational problem.
Let $\mathbb{P}_{M'}(\Lambda)$ denote the spaces of polynomials
of degree up to {$M'$} and set
$$ P_{M'}(\Lambda) = \{p(x) \in \mathbb{P}_{M'}(\Lambda)| p(0) = p(X) = 0\}.$$
Adopt the Fourier-like functions proposed by Zheng et al.\
\cite[see Section 4.2 for more details]{Zheng2015}.
\commenttwo{
Let $$\hat x = \frac{2}{X}(x - \frac{X}{2}),$$ and define 
$$z_k(x) = \lambda_k(L_k(\hat x) - L_{k+2}(\hat x)),$$
where $L_k(\hat x)$ is Legendre polynomials and $\lambda_k = \sqrt{\frac{X}{4(2k+3)}}$.
Then $(z'_i(x), z'_j(x))_\Lambda = \delta_{ij}$ where $\delta_{ij}$ is the Kronecker delta.
Set $Z = \left( (z_i, z_j)_\Lambda \right)_{(M'-1)\times(M'-1)}$,
and let $$\zeta_i(x) = \sum_{i=0}^{N-2}q_{ik} z_i(x), \quad k = 0, 1, \cdots, M'-2.$$
where $\{q_{ik}\}_{i=0}^{M'-2}$ is the eigenvector corresponding to the eigenvalue $\pi_k$ of the matrix $Z$.}
And the basis functions $\zeta_i(x)$ have the following property.
\begin{lemma}[\cite{Zheng2015}]
  For the basis functions $\zeta_i(x)(i=0,1,\cdots,M-2),x\in(a,b)$,
  \begin{equation}
    (\zeta_i(x), \zeta_j(x)) = \pi_i \delta_{ij}, \quad
    (\zeta_i'(x), \zeta_j'(x)) = \delta_{ij}.
  \end{equation}
\end{lemma}
Let $v_L^n = \sum_{i=0}^{M'-2}\hat v_i^n\zeta_i$, then we obtain the spectral scheme
for problem \eqref{eq:var}:
\begin{equation}\label{eq:sscheme}
  \begin{aligned}
    \big( v^n_L(x),\zeta_{i}(x)\big) +
    &\ K_\gamma\sum_{j=0}^{N}r_{nj}\big( v^j_{L,x}(x),\zeta_{i}'(x)\big)  \\
    &=\big(h^n(x),\zeta_{i}(x)\big) +  \sum_{j=0}^{N}r_{nj}\big(g^j(x),\zeta_{i}(x)\big), \\
   &\qquad\qquad\qquad\qquad\forall i=0,1,\cdots, M'-2.
  \end{aligned}
\end{equation}
Now, we introduce the matrix form of this scheme.
Adopt the following notations
\begin{equation}
  \begin{aligned}
    & V^j = (\hat v_0^j,\hat  v_1^j, \cdots,\hat  v_{M'-2}^j)^T,\\
    & \commentone{ H^j = \big((h^j, \zeta_0), (h^j, \zeta_1), \cdots, (h^j, \zeta_{M'-2})\big)^T}, \\
    & \commentone{ G^j = \big((g^j, \zeta_0), (g^j, \zeta_1), \cdots, (g^j, \zeta_{M'-2})\big)^T}. \\
  \end{aligned}
\end{equation}
Then we have
\commentone{%
\begin{equation}
  \bm\Pi V^n + K_\gamma\sum_{j=0}^{N} r_{nj} V^j =  H^n + \sum_{j=0}^{N} r_{nj} G^j,
  \quad n = 0, 1, \cdots, N,
\end{equation}%
}%
where $\bm\Pi = \text{diag}(\pi_0,\pi_1,\cdots, \pi_{M'-2})$.
Also define
\begin{equation}
  \begin{aligned}
    &V = \Big({V^0}^T, {V^1}^T, \cdots, {V^N}^T\Big)^T,\\
    &H = \Big({H^0}^T, {H^1}^T, \cdots, {H^N}^T\Big)^T,\\
    &G = \Big({G^0}^T, {G^1}^T, \cdots, {G^N}^T\Big)^T.\\
  \end{aligned}
\end{equation}
Then the matrix form of this scheme can be written as
\begin{equation}
  (\bm A + K_\gamma \bm B) V = H + \bm B G,
\end{equation}
where $\bm A = \bm I_{N+1}\otimes \bm\Pi $, $\bm B = \bm W\otimes \bm I_{M'-1}$,
and $\bm I_{N+1}$, along with $\bm I_{M'-1}$, are identity matrices.
\section{Numerical experiments}\label{sec:example}
\commentTwo{Aim to verify the validity of our schemes, several test problems
  are presented in this section. The first two, of which the exact solutions
  are known, are respectively adopted to illustrate the accuracy of
  scheme~\eqref{eq:compact} and scheme~\eqref{eq:sscheme}. A comparison between
scheme~\eqref{eq:sscheme} and algorithm of \cite{Zheng2015} is also given in
the second example.  And the last one with an unknown exact solution shows the
behaviors of the sub-diffusion system.}
\begin{example}\label{example:sin}
  Consider equation \eqref{eq:orig} with $X = 1, T = 1$, and
  \begin{equation}
    \left\{
    \begin{aligned}
      &u(x,0) = 0, \\
      &u(0,t)= 0,\quad u(1,t) = t^{c+\gamma}\sin1, \\
      &f(x,t) =\left(k_\gamma t^{c}+t^{c+\gamma }\right)\sin x,
    \end{aligned}
    \right.
  \end{equation}
  where $k_\gamma=\frac{\Gamma(c+\gamma+1)}{\Gamma(c+1)}$ and $K_\gamma=1$\ \commenttwo{\cite{Gao2015}}.
  The exact solution under these conditions is 
  $u(x,t) = t^{c+\gamma}\sin x$.
\begin{table}
  \caption{Errors and spatial convergence order of \commenttwo{scheme~\eqref{eq:compact}} with 
  $c = 1.9$, $\gamma = 0.6$, $N = 100$}\label{table:sin1so1}
  \centering
    \small
  \begin{tabular}{ccccccc}
    \toprule
       &\multicolumn{2}{c}{$q = 1$}&\multicolumn{2}{c}{$q = 2$}
       &\multicolumn{2}{c}{$q = 3$}\\ \cmidrule{2-7}
       $M$& error & order& error & order& error & order\\
    \midrule
          10 &   2.11e-08 &            &   2.11e-08 &            &   2.11e-08 &            \\
          20 &   1.32e-09 &       4.00 &   1.32e-09 &       4.00 &   1.32e-09 &       4.00 \\
          30 &   2.60e-10 &       4.00 &   2.60e-10 &       4.00 &   2.60e-10 &       4.00 \\
          40 &   8.23e-11 &       4.00 &   8.23e-11 &       4.00 &   8.23e-11 &       4.00 \\
          50 &   3.37e-11 &       4.00 &   3.37e-11 &       4.00 &   3.37e-11 &       4.00 \\
    \bottomrule
  \end{tabular}
\end{table}
\begin{table}
  \caption{Errors and spatial convergence order of \commenttwo{scheme~\eqref{eq:compact}} with 
  $c = 0.1$, $\gamma = 0.6$, $N = 100$}\label{table:sin1so2}
  \centering
    \small
  \begin{tabular}{ccccccc}
    \toprule
       &\multicolumn{2}{c}{$q = 1$}&\multicolumn{2}{c}{$q = 2$}
       &\multicolumn{2}{c}{$q = 3$}\\ \cmidrule{2-7}
       $M$& error & order& error & order& error & order\\
    \midrule
            10 &   2.58e-08 &            &   2.27e-08 &            &   2.27e-08 &            \\
            20 &   4.76e-09 &       2.44 &   1.42e-09 &       4.00 &   1.42e-09 &       4.00 \\
            30 &   4.21e-09 &       0.30 &   2.80e-10 &       4.00 &   2.81e-10 &       4.00 \\
            40 &   4.30e-09 &      -0.07 &   8.86e-11 &       4.00 &   8.88e-11 &       4.00 \\
            50 &   4.23e-09 &       0.07 &   3.62e-11 &       4.02 &   3.64e-11 &       4.00 \\
    \bottomrule
  \end{tabular}
\end{table}
  In this example, we investigate the convergence order of the
  scheme~\eqref{eq:compact}. The high convergence order can
  be easily seen from Table \ref{table:sin1so1}, in which we take $c= 1.9$,
  $\gamma=0.6$, and $N=100$.  Table \ref{table:sin1so2} shows the effect of the 
  smoothing method with different $q$ for equation~\eqref{eq:orig} with $c =
  0.1$,$\gamma=0.6$, and $N=100$.  It also shows that the convergence order in
  space is significantly improved when $q=2$.
  As shown in \commentthree{Tables~\ref{table:sin1to1} and \ref{table:sin1to2}}, the convergence
  rate in time can be improved in several degrees when the exact solutions have
  different regularity.
  These results indicate that the smoothing method can retain the
  convergence order when the regularity of $u$ is low.
\begin{table}
  \caption{Maximum errors of \commenttwo{scheme~\eqref{eq:compact}} with
  $c = 1.4$, $\gamma = 0.8$, $M = 5000$}\label{table:sin1to1}
  \centering
    \small
  \begin{tabular}{cccccc}
    \toprule
           $N$ &          6 &          8 &         10 &         12 &         14 \\
    \midrule
       $q = 1$ &   3.20e-06 &   6.91e-07 &   2.09e-07 &   7.76e-08 &   3.33e-08 \\
       $q = 2$ &   2.37e-07 &   1.08e-08 &   9.91e-10 &   1.42e-10 &   3.12e-11 \\
       $q = 3$ &   5.44e-07 &   3.26e-09 &   7.84e-11 &   2.13e-12 &   1.45e-12 \\
    \bottomrule
  \end{tabular}
\end{table}
\begin{table}
  \caption{Maximum errors of \commenttwo{scheme~\eqref{eq:compact}} with 
  $c = 0.5$, $\gamma = 0.8$, $M = 5000$}\label{table:sin1to2}
  \centering
    \small
  \begin{tabular}{cccccc}
    \toprule
           $N$ &          6 &          8 &         10 &         12 &         14 \\
    \midrule
       $q = 1$ &   2.87e-05 &   1.00e-05 &   4.39e-06 &   2.22e-06 &   1.24e-06 \\
       $q = 2$ &   1.68e-07 &   1.84e-08 &   3.45e-09 &   8.72e-10 &   2.70e-10 \\
       $q = 3$ &   3.96e-07 &   2.22e-08 &   2.19e-09 &   3.42e-10 &   6.34e-11 \\
    \bottomrule
  \end{tabular}
\end{table}
\begin{table}
  \caption{Maximum errors of \commenttwo{scheme~\eqref{eq:compact}} with 
  $c = 3.1$, $\gamma = 0.5$, $M = 5000$}\label{table:sin1to3}
  \centering
    \small
  \begin{tabular}{cccccc}
    \toprule
           $N$ &          6 &          8 &         10 &         12 &         14 \\
    \midrule
       $q = 1$ &   7.23e-08 &   6.42e-09 &   1.03e-09 &   2.82e-10 &   1.03e-10 \\
       $q = 2$ &   4.90e-06 &   5.38e-09 &   4.49e-11 &   4.82e-12 &   2.46e-12 \\
       $q = 3$ &   2.85e-04 &   4.79e-06 &   3.32e-09 &   1.64e-11 &   3.13e-12 \\
    \bottomrule
  \end{tabular}
\end{table}
\begin{figure}
    \centering
    \subfloat[$c=1.4$]{
    \includegraphics[width=0.45\textwidth]{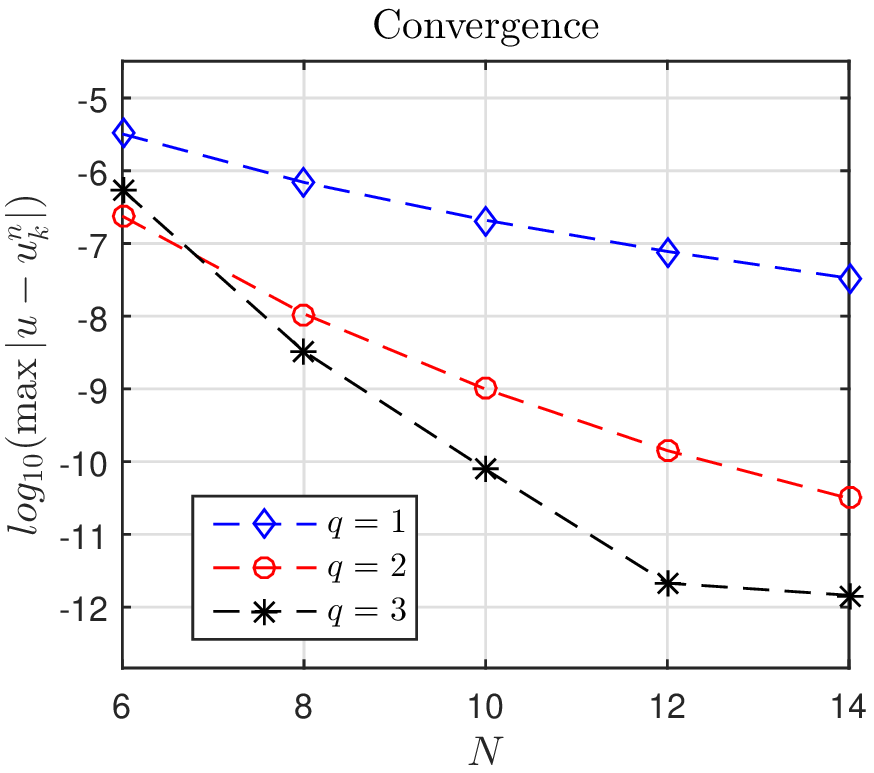}
  }
    \subfloat[$c=0.5$]{
    \includegraphics[width=0.45\textwidth]{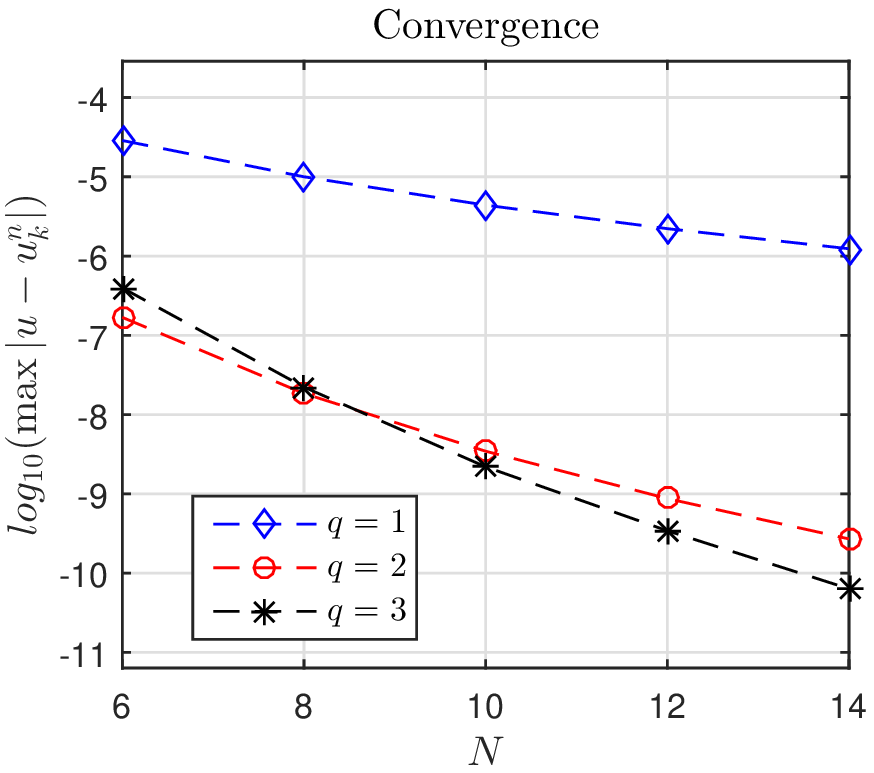}
  }
    \caption{Convergence of scheme~\eqref{eq:compact} for different $c$ with $\gamma = 0.8$, $M = 10000$.}\label{figure:sin1_1}
\end{figure}
\begin{figure}
    \centering
    \includegraphics[width=0.45\textwidth]{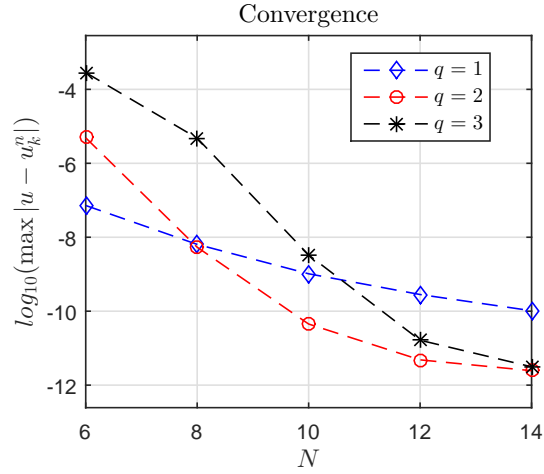}
    \caption{Convergence of scheme~\eqref{eq:compact}  with $c = 3.1$, $\gamma = 0.5$, $M = 10000$.}\label{figure:sin1_2}
\end{figure}
  Maximum errors for different $N$ are shown in Figure~\ref{figure:sin1_1}
  with $c = 1.4$, $\gamma=0.8$ and $c = 0.5, \gamma =0.8$ separately.
  The Figure~\ref{figure:sin1_2} displays that the convergence rate \commentthree{increases as $q$ increases}.
  The results shown in these figures, along with the results in Table~\ref{table:sin1to3},
  \commentthree{suggest} that it is profitable to perform smooth transformation on the equation
  no matter what regularity the solution has.
\end{example}
\begin{example}\label{example:compare}
  In this example, we
  \commentthree{consider} equation \eqref{eq:orig} with $X = 1, T = 1$, and
  \begin{equation}
    \left\{
    \begin{aligned}
      &u(x, 0)= 0, \\
      &u(0,t) = 0,\quad u(1,t) = 0, \\
      &f(x,t) =\left(k_\gamma t^{c}+\pi^2 t^{c+\gamma }\right)\sin \pi x,
    \end{aligned}
    \right.
  \end{equation}
  where $k_\gamma=\frac{\Gamma(c+\gamma+1)}{\Gamma(c+1)}$ and $K_\gamma=1$.
  The exact solution of equation \eqref{eq:orig} is
  $$u(x,t) = t^{c+\gamma}\sin \pi x.$$
  It is also the solution of equation in the following form:
  \begin{equation}\label{eq:zheng}
    \left\{
    \begin{aligned}
      &\prescript{RL}{0}{\mathcal D}^{\gamma}_T u = u_{xx} + f \\
      &u(x,0) = 0, u(t, 0) = 0, u(t, 1) = 0, \\
      &f = k_\gamma t^{c}\sin \pi x + \pi^2 t^{c+\gamma }\sin \pi x,
    \end{aligned}
    \right.
  \end{equation}
  This equation is equivalent with our equation \eqref{eq:orig}.
  In \cite{Zheng2015}, the authors solved these sub-diffusion equations with this form.
  Here, we compare the numerical results of scheme~\eqref{eq:sscheme} with
  the results of the method developed by Zheng et al.\ in \cite{Zheng2015},
  which possesses high efficiency and exponential decay in both time and space directions.
\begin{table}
  \caption{Maximum errors of \commenttwo{scheme~\eqref{eq:sscheme}} with 
  $c=2.5$, $\gamma = 0.4$, $M'=200$}\label{table:sp1to1}
  \centering
    \small
  \begin{tabular}{cccccc}
    \toprule
           $N$ &          6 &          8 &         10 &         12 &         14 \\
    \midrule
       $q = 1$ &   4.87e-07 &   6.38e-08 &   1.26e-08 &   3.23e-09 &   1.00e-09 \\
       $q = 2$ &   5.31e-07 &   7.01e-09 &   1.80e-10 &   7.27e-12 &   4.90e-13 \\
       $q = 3$ &   7.60e-05 &   8.67e-08 &   7.51e-10 &   1.56e-11 &   3.13e-13 \\
   Zheng et al.& 6.66e-06 &   9.98e-07 &   2.26e-07 &   6.62e-08 &   2.32e-08 \\
    \bottomrule
  \end{tabular}
\end{table}
\begin{table}
  \caption{Maximum errors of \commenttwo{scheme~\eqref{eq:sscheme}} with 
  $c = 1.5$, $\gamma = 0.4$, $M'=200$}\label{table:sp1to2}
  \centering
    \small
  \begin{tabular}{cccccc}
    \toprule
           $N$ &          6 &          8 &         10 &         12 &         14 \\
    \midrule
       $q = 1$ &   2.25e-06 &   5.10e-07 &   1.54e-07 &   5.58e-08 &   2.32e-08 \\
       $q = 2$ &   2.39e-07 &   7.95e-09 &   4.85e-10 &   5.31e-11 &   8.67e-12 \\
       $q = 3$ &   8.45e-07 &   1.28e-08 &   3.80e-10 &   2.61e-11 &   1.17e-12 \\
    Zheng et al.&3.60e-05 &   9.63e-06 &   3.38e-06 &   1.41e-06 &   6.67e-07 \\
    \bottomrule
  \end{tabular}
\end{table}
\begin{figure}
    \centering
    \subfloat [$c=2.5$]{
    \includegraphics[width=0.45\textwidth]{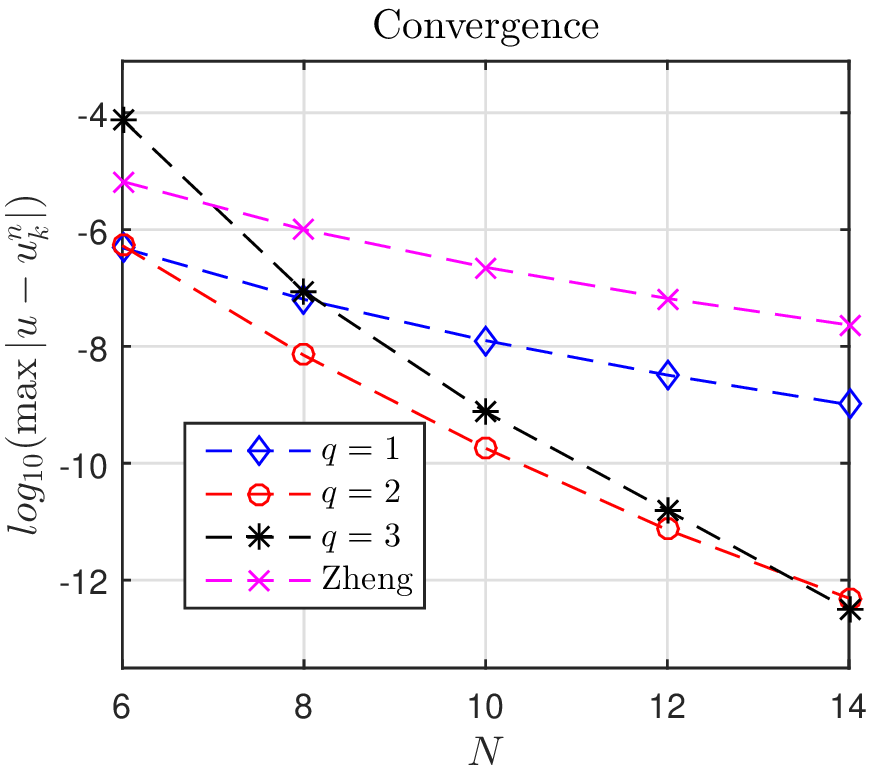}
  }
    \subfloat [$c=1.5$]{
    \includegraphics[width=0.45\textwidth]{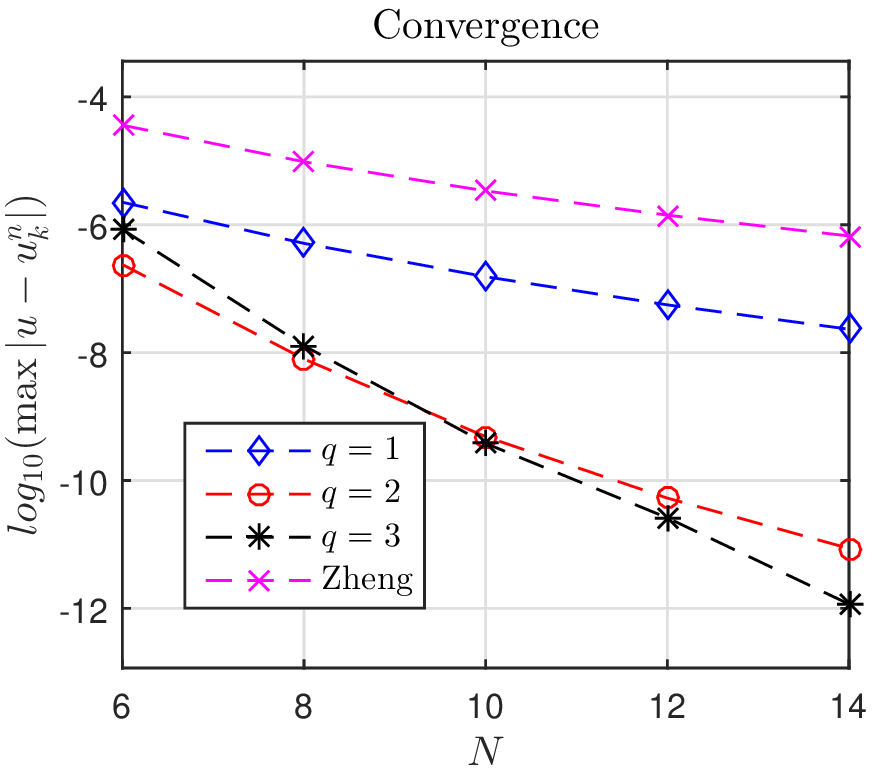}
  }
    \caption{Convergence of scheme~\eqref{eq:sscheme} for different $c$ with
    $\gamma = 0.4$, $M' = 200$.}\label{figure:sp1_1}
\end{figure}
  In this example we take $\gamma = 0.4$ and use the polynomials degree of
  spatial base $M' = 200$ to calculate results with different $c$.
  First,  we choose $c = 2.5, 1.5$, and present the maximum errors in
  Table~\ref{table:sp1to1} and Table~\ref{table:sp1to2} respectively.
  In addition to visualize these data, we plot them in Figure~\ref{figure:sp1_1}.
  From those tables and figures, we can see that Zheng's result is close to our
  result with $q = 1$, and larger $q$ can significantly improve the accuracy.
  Also, notice that larger $q$ can speed up the convergence.
\begin{table}
  \caption{Maximum errors of \commenttwo{scheme~\eqref{eq:sscheme}}  with 
  $c = 0.5$, $\gamma = 0.4$, $M'=200$}\label{table:sp1to3}
  \centering
    \small
  \begin{tabular}{cccccccc}
    \toprule
           $N$ &         10 &         18 &         26 &         34 &         42 &         50 \\
    \midrule
       $q = 1$ &   4.38e-06 &   4.64e-07 &   1.04e-07 &   3.39e-08 &   1.40e-08 &   6.76e-09 \\
       $q = 2$ &   1.27e-08 &   1.42e-10 &   7.00e-12 &   7.23e-13 &   1.19e-13 &   2.18e-14 \\
       $q = 3$ &   4.31e-09 &   1.51e-11 &   3.63e-13 &   2.32e-14 &   1.11e-15 &   1.22e-15 \\
     Zheng et al.&1.48e-04 &2.56e-05 &8.09e-06 &3.42e-06 &1.72e-06 &9.68e-07 \\
    \bottomrule
\end{tabular}
\end{table}
\begin{table}
  \centering
  \caption{Maximum errors of \commenttwo{scheme~\eqref{eq:sscheme}}  with 
  $c = -0.1$, $\gamma = 0.4$, $M'=200$}\label{table:sp1to4}
    \small
  \begin{tabular}{ccccccc}
    \toprule
           $N$ &         10 &         18 &         26 &         34 &         42 &         50 \\
    \midrule
       $q = 1$ &   1.95e-05 &   4.14e-06 &   1.50e-06 &   7.15e-07 &   4.05e-07 &   2.57e-07 \\
       $q = 2$ &   1.47e-06 &   9.93e-08 &   1.68e-08 &   4.50e-09 &   1.58e-09 &   6.71e-10 \\
       $q = 3$ &   1.57e-07 &   3.65e-09 &   3.22e-10 &   5.32e-11 &   1.27e-11 &   3.83e-12 \\
    \bottomrule
  \end{tabular}
\end{table}
\begin{figure}
    \centering
    \subfloat [$c=0.5$]{
    \includegraphics[width=0.45\textwidth]{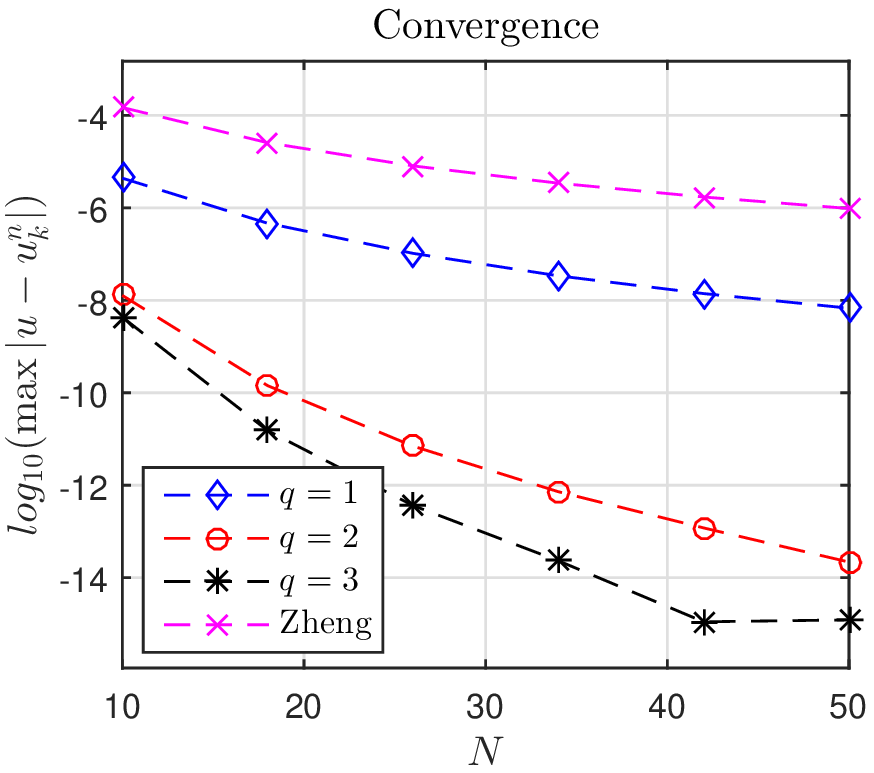}
  }
    \subfloat [$c=-0.1$]{
    \includegraphics[width=0.45\textwidth]{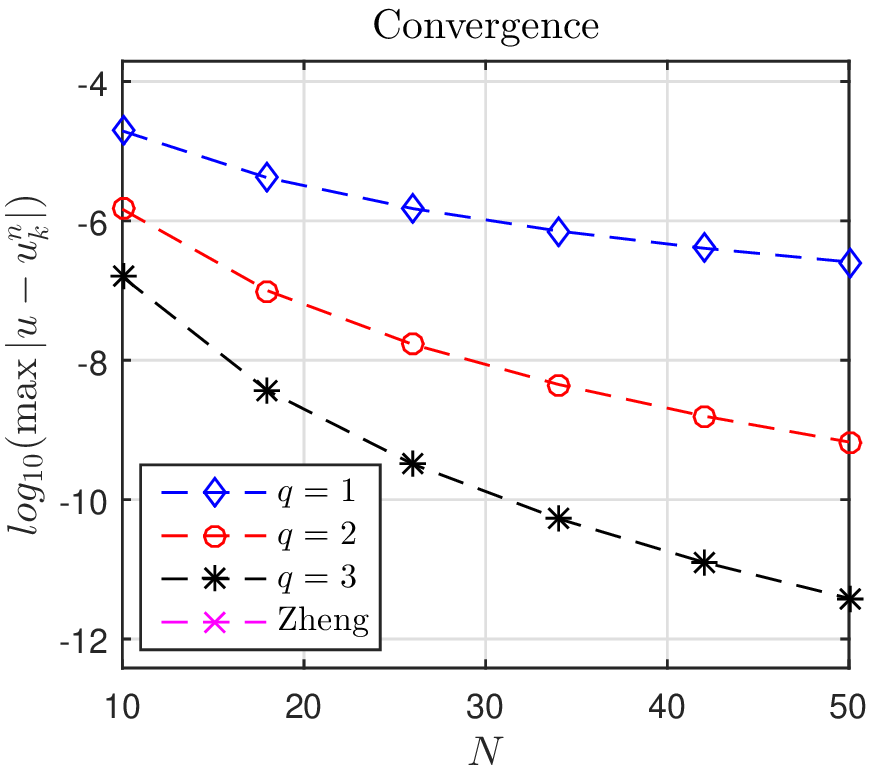}
  }
    \caption{Convergence of scheme~\eqref{eq:sscheme} for different $c$ with
    $\gamma = 0.4$, $M' = 200$.}\label{figure:sp1_2}
\end{figure}
  Furthermore, we choose $c = 0.5, -0.1$, and present the results in
  Table~\ref{table:sp1to3}, Table~\ref{table:sp1to4}, and Figure~\ref{figure:sp1_2}.
  Zheng's method cannot handle the situation with $c = -0.1$ and hence we just present
  our result in table and figure when $c = -0.1$.
  As $c$ \commentthree{decreases}, the regularity of the solution \commentthree{becomes} weaker, 
  and this \commentthree{results} in
  slow convergence speed as shown in Figure~\ref{figure:sp1_1} and Figure~\ref{figure:sp1_2}.
  However, we can still get accurate results by setting $q$ larger.  It is recommended
  to apply smooth transformation on the equations.
\end{example}
\begin{example}\label{example:real}
  Consider anomalous sub-diffusion equation
  \begin{equation}
    \frac{\partial u}{\partial t}=\frac{\partial^{1-\gamma}}{\partial t^{1-\gamma}}
    \Big(\frac{\partial^2u}{\partial x^2}\Big),
  \end{equation}
  with initial and boundary conditions \cite{Zhuang2008}
  \begin{equation}
    \begin{aligned}
      &u(x,0) = \left\{
        \begin{aligned}
          &2x, & 0 \le x \le 0.5, \\
          &\frac{4-2x}{3}, & 0.5\le x \le 2,
        \end{aligned}\right.\\
      &u(0,t) = u(2, t) = 0, \quad 0 \le t \le 0.4\,.
    \end{aligned}
  \end{equation}
  This system is simulated by applying scheme~\eqref{eq:compact} with $N=20, M =20$.
  Figure~\ref{fig:real1} shows the numerical approximation $u(x,t)$
  when $\gamma=0.1$ and $\gamma=0.9$ respectively. And Figure~\ref{fig:real3}
  illustrates the change of approximation $u(x,t)$ as $\gamma$ vary in quantity.
  These results \commentthree{show} that the system exhibits sub-diffusion
  behaviors and the solution continuously depends on the time fractional
  derivative.
\begin{figure}
    \centering
    \subfloat[$\gamma=0.1$]{
      \includegraphics[width=0.45\textwidth]{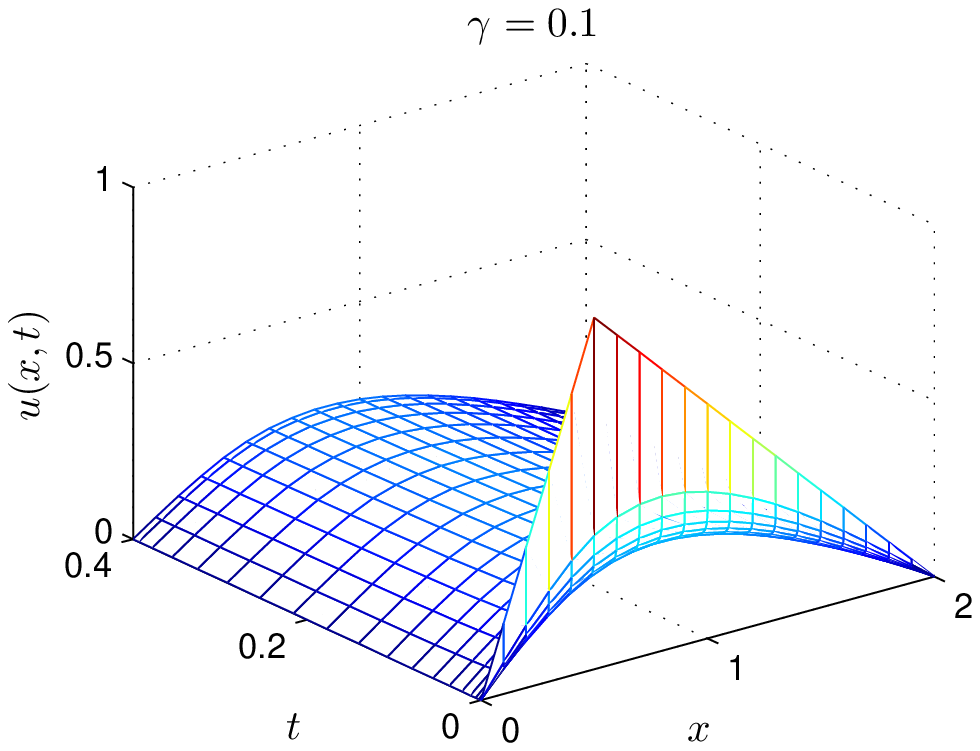}
    }
    \subfloat[$\gamma=0.9$]{
    \includegraphics[width=0.45\textwidth]{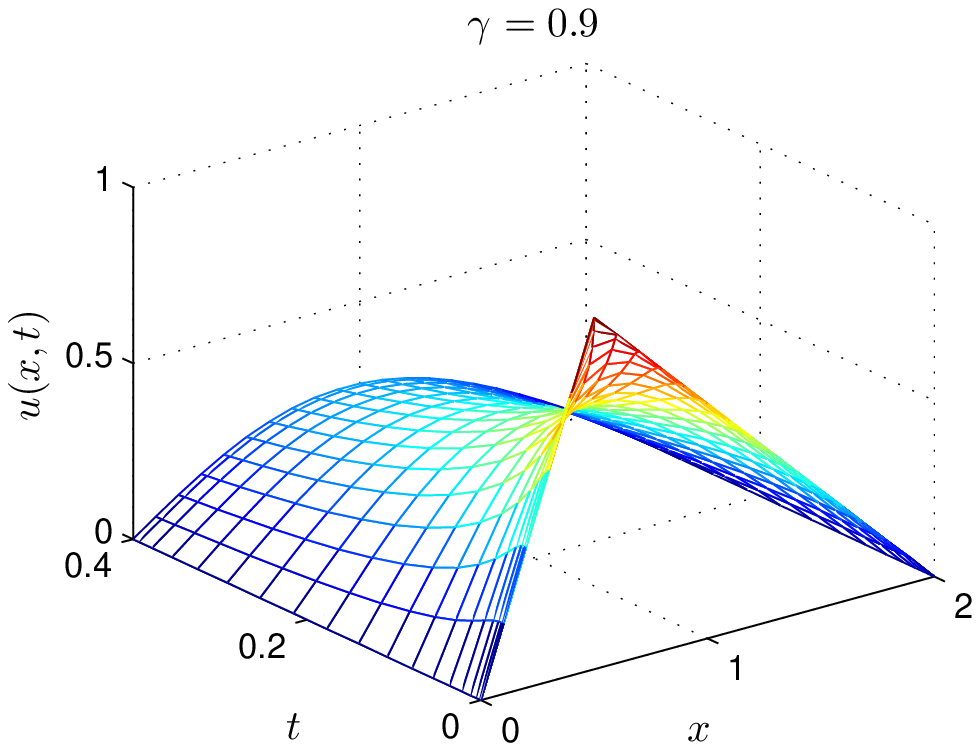}
    }
    \caption{Numerical approximation of $u(x,t)$ for different $\gamma$
     when $N = 20, M = 20$.}\label{fig:real1}
\end{figure}
\begin{figure}
    \centering
    \subfloat[$x=1.2$]{
      \includegraphics[width=0.45\textwidth]{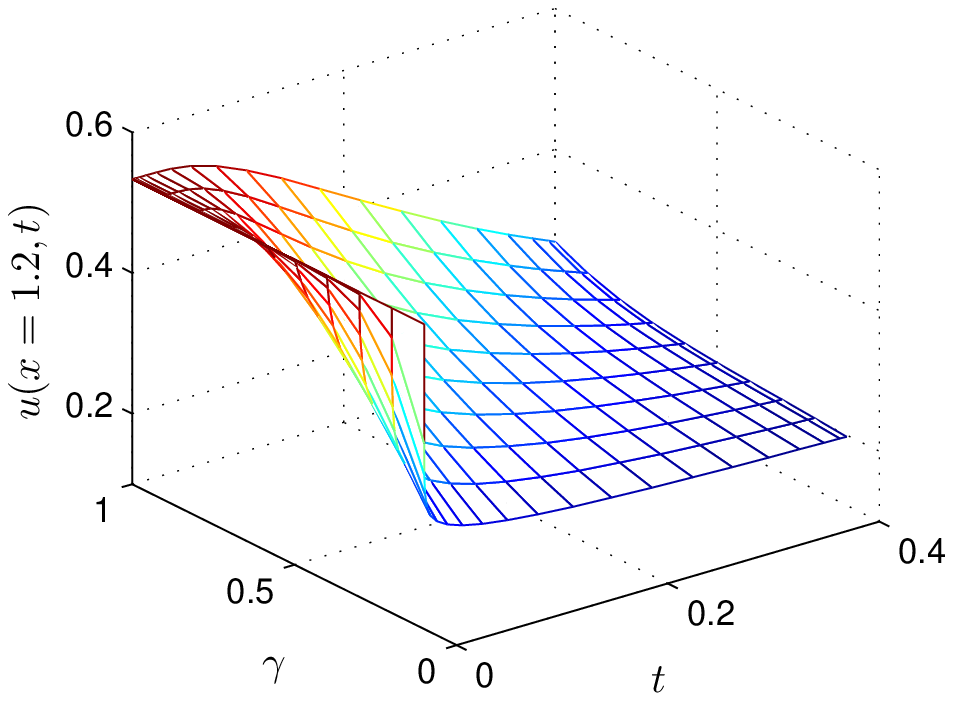}
    }
    \subfloat[$t=0.4$]{
      \includegraphics[width=0.45\textwidth]{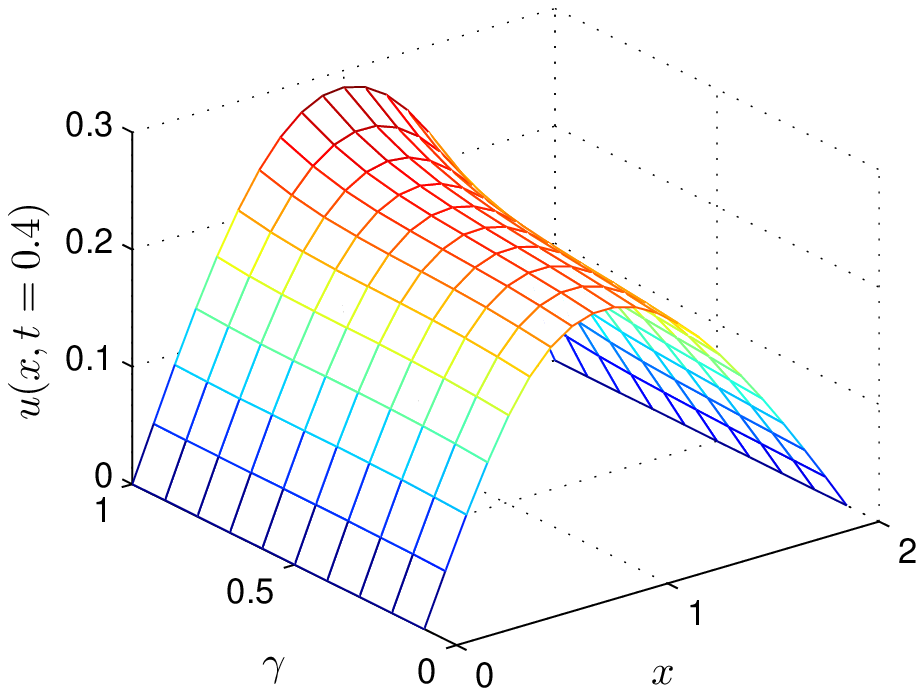}
    }
    \caption{The numerical approximation $u(x,t)$ for various $\gamma$ when $N = 20, M = 20$.}\label{fig:real3}
\end{figure}
\end{example}
\section{Conclusion}
In this paper, a high-order method has been proposed to solve anomalous
sub-diffusion equations especially when the exact solution has lower regularity. 
\commentOne{The compact difference method and the spectral method used in
  spatial direction make this method more effective.  In the numerical
  experiments, we have demonstrated the effectiveness and accuracy of these
  proposed schemes, even though theoretical analysis of convergence order and
  stability for this method are lacked because of the difficulty of coupled integral
  and differential operator.}
In the future, we will focus on the theoretical analysis of this method.
\section*{Acknowledgements}
The authors would like to thank the editor and the anonymous referees for their valuable comments and helpful suggestions that improve the quality of our paper. This research was supported by the National Natural Science Foundation of China (11601432,11471262) and the Strategic Research Grant of the City University of Hong Kong (7004446).
\bibliographystyle{plain}

\end{document}